\documentclass[twoside]{article}
\usepackage{ecj,palatino,epsfig,latexsym,natbib}

\usepackage{amsmath}
\usepackage{amssymb}
\usepackage{cleveref}
\usepackage{url}
\usepackage[dvipsnames,svgnames]{xcolor}
\usepackage{tikz}
\usetikzlibrary{arrows}
\usetikzlibrary{arrows.meta}
\usepackage{subcaption}

\usepackage{wasysym}

\tikzset{>={Triangle[scale=0.75]}}

\usepackage{algorithm}
\usepackage[noend]{algpseudocode}
\usepackage{algorithmicx}

\algblockdefx{Input}{EndInput}{\textbf{Input: }}{}
\algblockdefx{Output}{EndOutput}{\textbf{Output: }}{}
\algnotext{EndInput}
\algnotext{EndOutput}

\begin{document}

\ecjHeader{x}{x}{xxx-xxx}{202X}{R2 v2}{L.~Schäpermeier, P.~Kerschke}
% \title{\bf R2 Special Issue}  
% \title{\bf R2 v2: The Exactness Strikes Back}
% \title{\bf R2-D2: On the Exact Bi-objective R2 Indicator }
% \title{\bf R2 v2: The Pareto-compliant R2 Indicator for Bi-objective Benchmarking} 
\title{R2 v2:\\The Pareto-compliant R2 Indicator for Better Benchmarking in Bi-objective Optimization}

\author{\name{\bf Lennart Schäpermeier} \hfill \addr{lennart.schaepermeier@tu-dresden.de}\\ 
        \addr{Big Data Analytics in Transportation, TU Dresden, 
        01062 Dresden, Germany \\
        ScaDS.AI Dresden/Leipzig, 01062 Dresden, Germany
        }
\AND
       \name{\bf Pascal Kerschke} \hfill \addr{pascal.kerschke@tu-dresden.de}\\
        \addr{Big Data Analytics in Transportation, TU Dresden, 
        01062 Dresden, Germany \\
        ScaDS.AI Dresden/Leipzig, 01062 Dresden, Germany
        }
}

\maketitle

\begin{abstract}

In multi-objective optimization, set-based quality indicators are a cornerstone of benchmarking and performance assessment. 
They capture the quality of a set of trade-off solutions by reducing it to a scalar number.
One of the most commonly used set-based metrics is the R2 indicator, which describes the expected utility of a solution set to a decision-maker under a distribution of utility functions.
Typically, this indicator is applied by discretizing the latter distribution, yielding a weakly Pareto-compliant indicator. In consequence, adding a nondominated or dominating solution to a solution set \textit{may} -- but does not have to -- improve the indicator's value.

In this paper, we reinvestigate the R2 indicator under the premise that we have a continuous, uniform distribution of (Tchebycheff) utility functions.
We analyze its properties in detail, demonstrating that this continuous variant is indeed Pareto-compliant -- that is, any beneficial solution will improve the metric's value.
Additionally, we provide efficient computational procedures that (a) compute this metric for bi-objective problems in $\mathcal O (N \log N)$, and (b) can perform incremental updates to the indicator whenever solutions are added to (or removed from) the current set of solutions, without needing to recompute the indicator for the entire set.
As a result, this work contributes to the state-of-the-art Pareto-compliant unary performance metrics, such as the hypervolume indicator, offering an efficient and promising alternative.

\end{abstract}


\section{Introduction}

When optimizing any system, there is often not just one objective, but multiple criteria required to assess the quality of a solution.
Rather than aggregating these different optimization objectives into one, e.g., by means of a linear combination of the individual objectives, the domain of multi-objective (MO) optimization aims to find a set of (Pareto-)optimal trade-off solutions to present to a decision-maker \citep{miettinen1999nonlinear}.
However, to benchmark MO optimizers and facilitate algorithm design, parameter tuning, and automated algorithm selection, quantifying the quality of trade-off solutions in a unary set-based performance indicator is often necessary.

To be reasonably interpretable, it is recommended that a MO performance indicator fulfills a property that is known as Pareto compliance \citep{grimme2021peeking}.
More precisely, a set-based performance indicator is called Pareto-compliant, if and only if its indicator value for set A is better than for set B when set A dominates set B.
In addition, an indicator is called weakly Pareto-compliant if its value for set A is not worse than for set B.

Up to now, the hypervolume (HV) indicator, and variants thereof, are the only set-based performance indicators that are recognized as truly Pareto-comp\-liant \citep{zitzler1998multiobjective,beume2009,guerreiro2021}.
The HV indicator computes the $m$-dim\-en\-sio\-nal hypervolume dominated by the solution set w.r.t. to a user-specified anti-optimal reference point.

In contrast, there are multiple families of weakly Pareto-com\-pli\-ant indicators.
Exemplary and widely used representatives are the IGD+ indicator \citep{ishibuchi2015modified}, requiring an (approximated or known) reference Pareto front, or the R2 indicator \citep{hansen1998,brockhoff2012properties}, requiring an ideal/utopian reference point as well as a sample of aggregation (or: utility) functions.

The R2 indicator, in particular, may be an attractive choice as it requires an ideal rather than an anti-optimal reference point.
In many problems, the ideal point is easier to find, e.g., in multi-objective machine learning problems where optimal, but not always anti-optimal, values for loss functions are available.
Also, solutions that do not dominate a chosen reference point may not contribute to the dominated hypervolume, and a reference point far away from the Pareto front (PF) tends to put a high weight on solutions at the PF's boundary, which is also often undesirable.
Another benefit arises when constructing test problems from objectives with known single-objective optima, which lack a natural upper bound for the reference point or a clear ``region of interest''.

While the unary R2 indicator was initially defined as an integral over a continuum of utility functions by \cite{hansen1998}, it is usually only discussed and applied in an approximate manner for which the distribution of utility functions is discretized \citep{brockhoff2012properties,trautmann2013r2}. 
The latter comes with the benefit of being flexible regarding the involved utility functions.
It also provides a convenient way to compute the indicator as an average of multiple utility functions, however, sacrificing Pareto compliance in the process.

In this paper, we consider the most common R2 indicator definition with a uniform distribution of Tchebycheff utility functions \citep{hansen1998,brockhoff2012properties}. 
Our main contribution is the methodology to compute this indicator exactly for a set of solutions in the bi-objective case, thereby preserving its Pareto compliance.
We achieve this by shifting perspective away from averaging over predefined utility functions towards computing the R2 indicator contributions of each individual nondominated point, thereby eliminating weaknesses in the R2 indicator's properties.
Additionally, we demonstrate the exact R2 indicator values for individual points and linear Pareto fronts, as well as the approximate nature of the discretization-based approach commonly used so far.
This paper is an extended version of \cite{schaepermeier2024reinvestigating}, which additionally provides the following contributions:
\begin{itemize}
    \item a correction regarding the exact R2 indicator values of concave fronts,
    \item a refined derivation of the continuous R2 indicator,
    \item the means to compute the exclusive contribution of a solution from the non-dominated set,
    \item an incremental update variant suitable for benchmarking scenarios, and
    \item a comparison between the iterative R2 and hypervolume indicators.
\end{itemize}
In the meantime, \cite{jaszkiewicz2024exact} have also published a more detailed analysis of the exact calculation of the R2 indicator.
They show the Pareto compliance of the indicator for more than two objectives based on the quick hypervolume algorithm but do not discuss the derivation of exact indicator values, or benchmarking-critical aspects such as its incremental computation.

This paper is structured as follows: We introduce multi-objective optimization and set-based performance assessment in \Cref{sec:background}.
Then, in \Cref{sec:method}, we derive the methodology for computing exact R2 indicator values, first for a single solution, then for a set of solutions, and finally present an efficient method to compute the indicator incrementally for a history of evaluated points.
\Cref{sec:experimental} presents some exemplary results regarding characteristics of the exact R2 indicator, and \Cref{sec:conclusion} concludes the paper with an outlook on future research avenues.

\section{Background} \label{sec:background}

We begin by introducing some fundamental aspects of multi-objective optimization and dominance relationships of multi-objective solutions.
Then, we will cover core aspects of set-based performance assessment in multi-objective optimization and its best known representative, the hypervolume indicator.
Finally, we introduce the R2 indicator with its most important properties for the discretized and the exact variant.

\subsection{Multi-objective Optimization}

In multi-objective (MO) optimization, we aim to (w.l.o.g.) minimize multiple conflicting objectives.
Commonly, a MO optimization problem (MOP) with $m$ objectives is given by an objective function $F: \mathcal X \mapsto \mathbb R^m$ where $\mathcal X$ represents the decision space.
The individual objectives are denoted as $f_i: \mathcal X \mapsto \mathbb R, i=1,\dots,m$ in this work.
Further, we are primarily considering the bi-objective setting ($m=2$).
Depending on the particular problem and decision space, there may be further constraints on admissible solutions.

A particular challenge posed by MOPs pertains to solution comparison.
While in single-objective optimization, solutions can be compared directly (either they have identical objective values or one is better than another), such immediate comparisons are not possible for all solutions of a MOP.
To solve this, we need the concept of dominance.
A solution $x$ \textit{dominates} another solution $y$ ($x \prec y$), iff $f_i(x) \leq f_i(y)$ for all $i$ and $f_i(x) < f_i(y)$ for at least one $i$.
A solution $x$ \textit{strongly dominates} another solution $y$ if the stronger condition $f_i(x) < f_i(y)$ holds for all~$i$.
A solution that dominates, but does not strongly dominate, another solution is also called \textit{weakly dominant}.
Finally, two solutions can be incomparable, that is, \textit{mutually nondominated}, if either fulfills some objective better than the other.

Definitions of dominance can also be extended to sets of solutions.
A set of solutions $A$ (weakly) dominates another set $B$, if each member of $B$ is (weakly) dominated by a solution in $A$, written as $A \preceq B$ and $A \prec B$, respectively \citep{zitzler2003performance,grimme2021peeking}.

The set of all nondominated solutions $P = \{x \in \mathcal X ~ | ~ \nexists y \in \mathcal X: y \prec x\}$ is known as the Pareto set, and its image under $F$ is known as the Pareto front.
The Pareto set contains the optimal trade-off solutions regarding the objectives, and aiming to obtain a close approximation to it is the prevalent approach of solving MOPs when no further constraints or preferences on the objectives are known, i.e., under black-box assumptions.
Evolutionary algorithms are the most widespread approach for finding good Pareto set approximations in these conditions.

Finally, we call the vector of the optimal, individual function values \emph{ideal point}, and refer to the best vector dominated by all Pareto optimal points as \emph{nadir point}.
Often, before computing indicator values and if the ideal and nadir points are available, the region between them is normalized to the $[0,1]^m$ box in objective space as a normalization technique.

\subsection{Set-based Performance Assessment}

As the Pareto set generally contains more than one solution, set-based performance measures are the norm in assessing the overall quality of an archive of evaluated points.
This need to quantify Pareto set approximations has led to numerous performance measures being introduced.
For a recent survey on MO performance indicators, we refer to \cite{audet2021performance}.
The most prominent set-based performance measure is the dominated hypervolume (HV) indicator (or: $S$-metric) that measures the space dominated by the set of solutions w.r.t.~an anti-optimal reference point \citep{zitzler1998multiobjective,beume2009,guerreiro2021}.
An illustration of the HV indicator for two objectives is given in \Cref{fig:hv}.

An important property of a set-based performance measure $I: \mathbb R^m \rightarrow \mathbb R$ is \emph{Pareto compliance}: If $A \preceq B$ and $B \npreceq A$, then $I(A) < I(B)$ \citep{zitzler2003performance}.
That is, a performance measure should improve if new non-dominated or (weakly) dominating solutions are added to a set of solutions.
If only $I(A) \leq I(B)$ can be guaranteed under the same circumstances, $I$ is called \emph{weakly Pareto-compliant}.
Only the HV indicator and other indicators based on it are established to be Pareto-compliant \citep{guerreiro2021}.
The selection of weakly Pareto-compliant indicators is somewhat larger, including, for example, the (discretized) R2 \citep{brockhoff2012properties} and the IGD+ \citep{ishibuchi2015modified} measures.

The necessity of the anti-optimal reference point for computing the HV indicator can, however, be a hindrance to achieving Pareto compliance in practice.
Setting the reference point so far back that it is dominated by every feasible solution introduces a bias towards the edges of the PF, while a reference point close to the nadir point fails to consider all solutions outside of such a defined region of interest \citep[see, e.g.,][]{ishibuchi2018}, cf.~\Cref{fig:hv}.

\begin{figure}[tb]
    \centering
    \begin{tikzpicture}        
        \coordinate (ref) at (3.4,3.4);
        
        \coordinate (lu) at (0.7,3.4);
        \coordinate (ym1) at (0.7,2.7);
        \coordinate (yint1) at (1.6,2.7);
        \coordinate (y) at (1.6,1.7);
        \coordinate (yint2) at (2.6,1.7);
        \coordinate (yp1) at (2.6,0.7);
        \coordinate (rl) at (3.4,0.7);
        
        \draw[->,ultra thick] (0,0)--(4,0) node[right, Orange]{$f_1$};
        \draw[->,ultra thick] (0,0)--(0,4) node[above, NavyBlue]{$f_2$};
        \filldraw[ultra thick, LightGray] (lu) -- (ym1) -- (yint1) -- (y) -- (yint2) -- (yp1) -- (rl) -- (ref) -- (lu);
        \draw[ultra thick, Gray] (lu) -- (ym1) -- (yint1) -- (y) -- (yint2) -- (yp1) -- (rl);
        
        \draw (ym1) node {\huge +};
        \draw (y) node {\huge +};
        \draw (yp1) node {\huge +};
        \draw[Gray] (ref) node {\huge +};
    \end{tikzpicture}
    \quad
    \begin{tikzpicture}        
        \coordinate (ref) at (3.4,3.4);
        
        \coordinate (lu) at (0.4,3.4);
        \coordinate (rl) at (3.4,0.4);
        \coordinate (ll) at (0.4,0.4);
        
        \draw[->,ultra thick] (0,0)--(4,0) node[right, Orange]{$f_1$};
        \draw[->,ultra thick] (0,0)--(0,4) node[above, NavyBlue]{$f_2$};
        \draw[ultra thick, dotted, LightGray] (ref) -- (lu) -- (ll) -- (rl) -- (ref) node[Gray]{\huge +};
        \draw[ultra thick, dashed, LightGray] (1.2, 4.0) -- (1.2, 3.4) node[at start, Black]{\huge +};
        \draw[ultra thick, dashed, LightGray] (2.2, 3.9) -- (2.2, 3.4) node[at start, Black]{\huge +};
        \draw[ultra thick, dashed, LightGray] (3.8, 2.2) -- (3.4, 2.2) node[at start, Black]{\huge +};
        \draw[ultra thick, dashed, LightGray] (4.0, 0.7) -- (3.4, 0.7) node[at start, Black]{\huge +};
        % \draw[ultra thick, Gray] (lu) -- (ym1) -- (yint1) -- (y) -- (yint2) -- (yp1) -- (rl);

        % \draw (ym1) node {\huge +};
        % \draw (y) node {\huge +};
        % \draw (yp1) node {\huge +};
        % \draw[Gray] (ref) node {\huge +};
    \end{tikzpicture}
    
    \caption{Left: Illustration of the HV indicator. Right: If no point dominating the HV is found, the minimal distance to the region of interest (dotted) can be used as an additional indicator. The combined indicator is, however, not Pareto-compliant anymore.}
    \label{fig:hv}
\end{figure}

\subsection{The R2 Indicator}

In contrast to the HV indicator, the unary R2 indicator is ordinarily defined as the expected utility of the point set w.r.t. a distribution of utility functions $U$ \citep{hansen1998}.
In the most general case, for a set of solutions $Y$ of a MOP, we can define it as follows:
$$
R2(Y) := \int_{u \in U} \min_{y \in Y} u(y) \, du.
$$
The most common choice of a utility function is a Tchebycheff aggregation, which allows to reach all Pareto-optimal points depending on the chosen parametrization.
For a weight vector $w \in [0,1]^m$ with $\sum_{i=1}^m w_i = 1$ and a utopian vector $y^*$, it is given by
\begin{align*}
    u_w(y) &= \max_{i = 1,\dots,m}{w_i (y_i - y_i^*)} \\
           &= \max_{i = 1,\dots,m}{w_i y'_i}
\end{align*}
using $y'_i = y_i - y_i^*$ to shift the utopian point w.l.o.g. to the origin.
The distribution of utility functions is then usually chosen as uniform on the weight simplex.

In the bi-objective case, where $w_2 = 1 - w_1$ holds, this yields the following formula for R2:
\begin{align*}
    R2(Y) &= \int_0^1 \min_{y \in Y} u_w(Y) \, dw \\
          &= \int_0^1 \min_{y \in Y} \{\max(w y'_1, (1 - w) y'_2)\} \, dw.
\end{align*}
As there is no apparent way to calculate this property directly, it is generally approximated in a discrete manner by discretizing $U$ using $n = |W|$ weight vectors $w \in W$:
\begin{align*}
    R2(Y) \approx \frac{1}{|W|} \sum_{w \in W} \min_{y \in Y} u_w(y).
\end{align*}
As an example, the uniform weight distribution with size $n$ for the bi-objective case is given by \citep{brockhoff2012properties}
$$
W = \left\{(0, 1), \left(\frac{1}{n-1}, \frac{n-2}{n-1}\right), \dots, \left(\frac{n-2}{n-1},\frac{1}{n-1}\right), \left(1, 0\right)\right\}.
$$

The discretization is simultaneously a blessing and a curse: On the one hand, it provides an effective method of approximating the underlying exact R2 value with high precision.
On the other hand, this weakens the indicators' properties.
To optimize the discretized R2 indicator, one can consider at most $|W|$ points on the Pareto front, which optimize $u_w$ for each $w \in W$, respectively.
Additional nondominated solutions cannot contribute to the indicator value.
Further, an individual utility function $u_w$ may be optimized by a point that is only weakly Pareto optimal, but has the same utility (for this set of weights) as another, Pareto optimal point.
See, for example, the top left image of \Cref{fig:r2-levelsets}, where each point along the vertical lines would have identical utility values.
This places the discretized R2 indicator among the weakly Pareto-compliant indicators.

\begin{figure}[tb]
    \centering
    \begin{tikzpicture}        
        \coordinate (y) at (1.7,1.3);
        \coordinate (O) at (0.5,0.5);
        
        \coordinate (w) at (0.5,1.5);
        \draw[densely dashed,very thick,LightGray] (0.5,1.5)--(1.5,0.5);

        \draw[very thick,Gray] (0.5,0.5)--(0.5,3);
        \draw[very thick,Gray] (1.1,0.5)--(1.1,3);
        \draw[very thick,Gray] (1.7,0.5)--(1.7,3);
        
        \draw[->,very thick] (0,0)--(3,0) node[right, Orange]{$f_1$};
        \draw[->,very thick] (0,0)--(0,3) node[above, NavyBlue]{$f_2$};
        \draw[very thick,Orange] (y)--(1.7,3);
        \draw[very thick,NavyBlue] (y)--(3,1.3);

        \draw[-{Stealth},very thick,LightGray,densely dashed] (O)--(w) node[below left,Gray]{$w$};

        \draw (y) node {\huge +} node[above right]{$y$};
        \draw[Gray] (O) node {\huge +};
    \end{tikzpicture}
    \begin{tikzpicture}        
        \coordinate (y) at (1.7,1.3);
        \coordinate (O) at (0.5,0.5);
        
        \coordinate (w) at (0.9,1.1);
        \draw[densely dashed,very thick,LightGray] (0.5,1.5)--(1.5,0.5);

        \draw[very thick,Gray] (0.5,1.1) -- (w) -- (0.9,0.5);
        \draw[very thick,Gray] (0.5,1.1 + 0.6) -- (0.9 + 0.4,1.1 + 0.6) -- (0.9 + 0.4,0.5);
        \draw[very thick,Gray] (0.5,1.1 + 2 * 0.6) -- (0.9 + 2 * 0.4,1.1 + 2 * 0.6) -- (0.9 + 2 * 0.4, 0.5);
        
        \draw[->,very thick] (0,0)--(3,0) node[right, Orange]{$f_1$};
        \draw[->,very thick] (0,0)--(0,3) node[above, NavyBlue]{$f_2$};
        \draw[very thick,Orange] (y)--(1.7,3);
        \draw[very thick,NavyBlue] (y)--(3,1.3);

        \draw[-{Stealth[scale=0.75]},very thick,LightGray,densely dashed] (O)--(w) node[below left,Gray]{$w$};

        \draw (y) node {\huge +} node[above right]{$y$};
        \draw[Gray] (O) node {\huge +};
    \end{tikzpicture}
    \begin{tikzpicture}        
        \coordinate (y) at (1.7,1.3);
        \coordinate (O) at (0.5,0.5);
        
        \coordinate (w) at (0.5 + 0.6, 0.5 + 0.4);
        \draw[densely dashed,very thick,LightGray] (0.5,1.5)--(1.5,0.5);

        \draw[very thick,Gray] (0.5,0.5 + 0.4) -- (w) -- (0.5 + 0.6,0.5);
        \draw[very thick,Gray] (0.5,0.5 + 0.6) -- (0.5 + 0.9,0.5 + 0.6) -- (0.5 + 0.9,0.5);
        \draw[very thick,Gray] (0.5,0.5 + 0.8) -- (0.5 + 1.2,0.5 + 0.8) -- (0.5 + 1.2, 0.5);
        
        \draw[->,very thick] (0,0)--(3,0) node[right, Orange]{$f_1$};
        \draw[->,very thick] (0,0)--(0,3) node[above, NavyBlue]{$f_2$};
        \draw[very thick,Orange] (y)--(1.7,3);
        \draw[very thick,NavyBlue] (y)--(3,1.3);

        \draw[-{Stealth[scale=0.75]},very thick,LightGray,densely dashed] (O)--(w) node[below left,Gray]{$w$};

        \draw (y) node {\huge +} node[above right]{$y$};
        \draw[Gray] (O) node {\huge +};
    \end{tikzpicture}
    \begin{tikzpicture}        
        \coordinate (y) at (1.7,1.3);
        \coordinate (O) at (0.5,0.5);
        
        \coordinate (w) at (0.5 + 0.7, 0.5 + 0.3);
        \draw[densely dashed,very thick,LightGray] (0.5,1.5)--(1.5,0.5);

        \draw[very thick,Gray] (0.5,0.5 + 0.3) -- (w) -- (0.5 + 0.7,0.5);
        \draw[very thick,Gray] (0.5,0.5 + 0.55) -- (0.5 + 1.2833,0.5 + 0.55) -- (0.5 + 1.2833,0.5);
        \draw[very thick,Gray] (0.5,0.5 + 0.8) -- (0.5 + 1.866,0.5 + 0.8) -- (0.5 + 1.866,0.5);
        
        \draw[->,very thick] (0,0)--(3,0) node[right, Orange]{$f_1$};
        \draw[->,very thick] (0,0)--(0,3) node[above, NavyBlue]{$f_2$};
        \draw[very thick,Orange] (y)--(1.7,3);
        \draw[very thick,NavyBlue] (y)--(3,1.3);

        \draw[-{Stealth[scale=0.75]},very thick,LightGray,densely dashed] (O)--(w) node[below left,Gray]{$w$};

        \draw (y) node {\huge +} node[above right]{$y$};
        \draw[Gray] (O) node {\huge +};
    \end{tikzpicture}
    \begin{tikzpicture}        
        \coordinate (y) at (1.7,1.3);
        \coordinate (O) at (0.5,0.5);
        
        \coordinate (w) at (1.5,0.5);
        \draw[densely dashed,very thick,LightGray] (0.5,1.5)--(1.5,0.5);

        \draw[very thick,Gray] (0.5,0.5)--(3,0.5);
        \draw[very thick,Gray] (0.5,0.9)--(3,0.9);
        \draw[very thick,Gray] (0.5,1.3)--(3,1.3);
        
        \draw[->,very thick] (0,0)--(3,0) node[right, Orange]{$f_1$};
        \draw[->,very thick] (0,0)--(0,3) node[above, NavyBlue]{$f_2$};
        \draw[very thick,Orange] (y)--(1.7,3);
        \draw[very thick,NavyBlue] (y)--(3,1.3);

        \draw[-{Stealth},very thick,LightGray,densely dashed] (O)--(w) node[below left,Gray]{$w$};

        \draw (y) node {\huge +} node[above right]{$y$};
        \draw[Gray] (O) node {\huge +};
    \end{tikzpicture}
    
    \caption{
    Illustration of level sets of the Tchebycheff utility for five different weight vectors $w$.
    The utility value $u_w(y)$ is determined by the surface that the level set touches first: At vertical surfaces (see left and center image in the top row), the $u_w(y)$ is determined by the $f_1$ value while at horizontal surfaces (see bottom row) $u_w(y)$ depends on the $f_2$ value of $y$.
    At $y$ (the top right figure) both $w_1 y_1$ and $w_2 y_2$ are identical.
    Note: Weight vectors are illustrated to point towards their equilibrium between both objectives, i.e., $w_1 f_1 = w_2 f_2$.
    }
    \label{fig:r2-levelsets}
\end{figure}

A final property of the R2 indicator is that each solution $y$ from a nondominated set of solutions $Y$ is optimal in terms of utility compared to all other solutions from $Y$ for a particular weight $(w_1^*, 1 - w_1^*)$ such that \citep{brockhoff2012properties}
\begin{align*}
    w_1^* y'_1 &= (1 - w_1^*) y'_2 \\
    \Rightarrow w_1^* &= \frac {y'_2} {y'_1 + y'_2}.
\end{align*}
Consequently, we also have $w_2^* = \frac {y'_1} {y'_1 + y'_2}$.
Discussions around the R2 indicator and its application focus (almost) exclusively on its discrete variant \citep{brockhoff2012properties,zitzler2008} rather than on the original continuous definition of R2 \citep{hansen1998}.
The remainder of this paper is dedicated to a better understanding of this original definition, analyzing its properties, and detailing methods for computing it.

\section{The R2 Indicator for Continuous Utility Distributions} \label{sec:method}

In this section, we derive how the R2 indicator can be computed under the assumption of a continuous distribution of Tchebycheff utility functions for bi-objective problems.
We will start by analyzing the utility of an axis-parallel segment and a single solution point before extending the analysis to sets of solutions, including the exclusive contribution of individual points and an incremental procedure to efficiently compute the indicator for a history of evaluated points.
Further, we derive the computational complexity of the presented approaches.

\subsection{Tchebycheff Utility of Axis-parallel Pareto Front Segments}

\begin{figure}[tb]
    \centering

    \begin{minipage}[t]{0.475\textwidth}
        \centering
        \begin{tikzpicture}
            \coordinate (y) at (2.0,1.5);
            \coordinate (y1) at (2.0,0.0);
            \coordinate (y2) at (0.0, 1.5);
            \coordinate (y2p) at (0.0, 3.0);
            \coordinate (yp) at (2.0,3.0);
            \coordinate (O) at (0.5,0.5);
            
            \draw[->,ultra thick] (0,0)--(4,0) node[right, Orange]{$f_1$};
            \draw[->,ultra thick] (0,0)--(0,4) node[above, NavyBlue]{$f_2$};
            \draw[ultra thick] (y1)--(2.0,-0.1) node[below]{$y_1$};
            \draw[ultra thick] (y2)--(-0.1,1.5) node[left]{$y_2$};
            \draw[ultra thick] (y2p)--(-0.1,3.0) node[left]{$y'_2$};
            \draw[ultra thick, LightGray, loosely dotted] (y)--(y1);
            \draw[ultra thick, LightGray, loosely dotted] (y)--(y2);
            \draw[ultra thick, LightGray, loosely dotted] (yp)--(y2p);
            \draw[ultra thick, Orange] (yp) -- (y);
    
            \draw (y) node {\huge +} node[below left]{$y$};
            \draw (yp) node {\huge +} node[below left]{$y'$};
            \draw[Gray] (O) node {\huge +};
        \end{tikzpicture}
        \caption{An axis-parallel (vertical) line segment between $y$ and $y'$ with utility $\color{Orange} u(y_1, [y_2, y_2'])$. For horizontal segments, the computation is analogous.}
        \label{fig:r2-segment}
    \end{minipage}
    \hfill
    \begin{minipage}[t]{0.475\textwidth}
        \centering
        \begin{tikzpicture}        
            \coordinate (y) at (2.0,1.5);
            \coordinate (y1) at (2.0,0.0);
            \coordinate (y2) at (0.0, 1.5);
            \coordinate (lu) at (2.0,4.0);
            \coordinate (rl) at (4.0,1.5);
            \coordinate (O) at (0.5,0.5);
            
            \draw[->,ultra thick] (0,0)--(4,0) node[right, Orange]{$f_1$};
            \draw[->,ultra thick] (0,0)--(0,4) node[above, NavyBlue]{$f_2$};
            \draw[ultra thick] (y1)--(2.0,-0.1) node[below]{$y_1$};
            \draw[ultra thick] (y2)--(-0.1,1.5) node[left]{$y_2$};
            \draw[ultra thick, Orange] (lu) -- (y);
            \draw[ultra thick, NavyBlue] (y) -- (rl);
            
            \draw[ultra thick, LightGray, loosely dotted] (y)--(y1);
            \draw[ultra thick, LightGray, loosely dotted] (y)--(y2);
    
            \draw (y) node {\huge +} node[below left]{$y$};
            \draw[Gray] (O) node {\huge +};
        \end{tikzpicture}
        \caption{Special case of a Pareto front consisting of a single point $y$. The corresponding R2 utility thus is $R2(\{y\}) = \color{Orange}u(y_1, [y_2, \infty))\color{black} + \color{NavyBlue}u(y_2, [y_1, \infty))$.}
        \label{fig:r2-single}
    \end{minipage}
\end{figure}

We begin by demonstrating how we can compute the partial R2 utility for an axis-parallel segment of the PF, cf.~\Cref{fig:r2-segment}.
This will be the building block for deriving the R2 value for a single solution as well as a set of solutions.

Without loss of generality, let $y^* = (0,0)$ be the utopian point.
Let $[y_2,y'_2]$ with $0 \leq y_2 < y'_2$ and $y_1 \geq 0$ be an axis-parallel segment of the Pareto front.
Let $w_1 = \frac {y_2} {y_2 + y_1}$ and $w'_1 = \frac {y'_2} {y'_2 + y_1}$ be the weights for which the respective end points of the segment are optimal w.r.t.~the R2 utility.
As $y_2 < y'_2$, we also have $w_1 < w'_1$.
The segment then contributes the partial R2 utility for the weight range $[w_1,w'_1]$ w.r.t.~$y_1$:
\begin{align*}
    u(y_1, [y_2, y'_2])
      &= \int_{w_1}^{w'_1} w y_1 dw \\
      % &= y_1 \left[\frac 1 2 w^2\right]_{w_1}^{w'_1} \\
      &= \frac 1 2 y_1 ({w'_1}^2 - {w_1}^2) \\
      &= \frac 1 2 y_1 \left(\left(\frac {y'_2} {y_1 + y'_2}\right)^2 - \left(\frac {y_2} {y_1 + y_2}\right)^2\right).
\end{align*}
If we are in the extreme regions of the PF where $y'_2 \rightarrow \infty$, the integral and computation are still well-defined as $\lim_{y'_2 \rightarrow \infty} {w'_1} = 1$.

The computation for an axis-parallel segment $[y_1,y'_1]$ with $0 \leq y_1 < y'_1$ and $y_2 \geq 0$ is completely analogous, so that we can write:
\begin{align*}
    u(y_2, [y_1, y'_1]) = \frac 1 2 y_2 \left(\left(\frac {y'_1} {y'_1 + y_2}\right)^2 - \left(\frac {y_1} {y_1 + y_2}\right)^2\right).
\end{align*}
%
% Finally, we observe that for a given interval $I = [y_2,y'_2]$ and given $y_1$, and $\varepsilon > 0$, $u(y_1, I) > u(y_1 + \varepsilon, I)$.
% That is, shifting a segment away from the ideal point always decreases its corresponding utility.
% In the limit, this utility equals zero, i.e., $\lim_{\varepsilon \rightarrow \infty} u(y_1 + \varepsilon, I) = 0$.
Finally, we observe that for a given interval $I = [y_2,y'_2]$, $\lim_{y_1 \rightarrow \infty} u(y_1, I) = 0$, i.e., we can use $u(\infty, I) = 0$.

\subsection{R2 Utility for a Single Solution}

Now, let $Y = \{y\}$ be the set containing only one solution $y = (y_1, y_2)$ that is dominated by $y^* = (0,0)$, cf.~\Cref{fig:r2-single}.
The Pareto front corresponding to $Y$ then only contains two axis-parallel segments with one unbounded end each, and the corresponding utility can be computed as follows:
\begin{align*}
    R2(\{y\}) &= \int_0^1 \max(y_1 w, y_2 (1 - w)) \, dw \\
          &= u(y_1, [y_2,\infty)) + u(y_2, [y_1,\infty)) \\
          &= \int_{\frac{y_2}{y_1 + y_2}}^{1} w y_1 dw + \int_{\frac{y_1}{y_1 + y_2}}^{1} w y_2 dw \\
          &= \frac 1 2 y_1 \left(1 - \left(\frac {y_2} {y_1 + y_2} \right)^2 \right) + \frac 1 2 y_2 \left(1 - \left(\frac {y_1} {y_1 + y_2} \right)^2 \right)
\end{align*}
This simple case demonstrates that computing the exact R2 indicator for a set of solutions consists of the following steps:
First, we need to identify areas in which the utility value does not vary w.r.t. $y$ and is only sensitive to $w$.
Then, we compute the contribution of each of these areas to the final R2 value using the corresponding integral and finally sum everything up.

We can build on these observations to derive a general procedure to compute $R2(Y)$ for arbitrary solution sets $Y$.

\subsection{R2 for a Set of Solutions}

\begin{figure}[tb]
    \centering
    \begin{tikzpicture}        
        \coordinate (lu) at (1.3,4.3);
        \coordinate (ym1) at (1.3,3.3);
        \coordinate (yint1) at (2.2,3.3);
        \coordinate (y) at (2.2,2.0);
        \coordinate (yint2) at (3.2,2.0);
        \coordinate (yp1) at (3.2,1.0);
        \coordinate (rl) at (4.2,1.0);
        \coordinate (O) at (0.5,0.5);
        
        \draw[->,ultra thick] (0,0)--(4,0) node[right, Orange]{$f_1$};
        \draw[->,ultra thick] (0,0)--(0,4) node[above, NavyBlue]{$f_2$};
        \draw[ultra thick, dashed, LightGray] (lu) -- (ym1);
        \draw[ultra thick, LightGray] (ym1) -- (yint1);
        \draw[ultra thick, Orange] (yint1) -- (y);
        \draw[ultra thick, NavyBlue] (y) -- (yint2);
        \draw[ultra thick, LightGray] (yint2) -- (yp1);
        \draw[ultra thick, dashed, LightGray] (yp1) -- (rl);

        \draw (ym1) node {\huge +} node[below left]{$y^{(n-1)}$};
        % \draw (yint1) node {\huge +};
        \draw (y) node {\huge +} node[below left]{$y^{(n)}$};
        % \draw (yint2) node {\huge +};
        \draw (yp1) node {\huge +} node[below left]{$y^{(n+1)}$};
        \draw[Gray] (O) node {\huge +};

        % \draw[very thick, lightgray, -{Stealth}, dashed] (O) -- (yint1);
        % \draw[very thick, lightgray, -{Stealth}, dashed] (O) -- (y);
        % \draw[very thick, lightgray, -{Stealth}, dashed] (O) -- (yint2);
    \end{tikzpicture}
    
    % \caption{Schematic illustration of the integration ranges of a solution $y^{(n)}$. $y^{(n-1)}$, $y^{(n)}$, and $y^{(n+1)}$ are consecutive points in the solution set. $y^{(n-)}$ and $y^{(n+)}$ correspond to the corners in the PF that are adjacent to $y^{(n)}$. Their corresponding weight vectors $w^{(n-)}$ and $w^{(n+)}$ indicate the boundaries in which $y^{(n)}$ locally determines the PF. At $w^{(n)}$, the objective switches between the $f_1$- (vertical PF segment) and $f_2$-values (horizontal PF segment) of solution $y^{(n)}$.}
    \caption{Integration ranges surrounding a solution $y^{(n)}$. $y^{(n-1)}$, $y^{(n)}$, and $y^{(n+1)}$ are consecutive points in the solution set. The utility of the vertical (orange) segment is determined by the $\color{Orange} f_1$ value of $y^{(n)}$, while the utility of points on the horizontal (blue) segment is dependent on its $\color{NavyBlue} f_2$ value. The length of the segments depends on the neighbors of $y^{(n)}$ in the solution set.}
    \label{fig:r2-comp-helper}
\end{figure}

Let us now consider what happens when our solution set contains $N > 1$ solutions.
Let $Y = \{y^{(1)}, \dots, y^{(N)}\}$ be the set of \emph{nondominated} solutions ordered by ascending $y_1$ value.
In order to compute $R2(Y)$, we can consider the utilities from the respective axis-parallel segments that make up the Pareto front.
We can identify all segments by considering each solution $y^{(n)}$ and its respective neighbors $y^{(n-1)}$ and $y^{(n+1)}$, cf.~\Cref{fig:r2-comp-helper} for a visual aid.

For the {\color{Orange} vertical segment} neighboring $y^{(n)}$, we get the utility $\color{Orange} u(y_1^{(n)}, [y_2^{(n)}, y_2^{(n-1)}])$, while for the {\color{NavyBlue} horizontal segment}, we get $\color{NavyBlue} u(y_2^{(n)}, [y_1^{(n)}, y_1^{(n+1)}])$.
For the special cases $n=1$ or $n=N$ at the extreme ends of the Pareto front, where there is no neighbor in one direction, we can use $y_1^{(N+1)} \rightarrow \infty$ and $y_2^{(0)} \rightarrow \infty$, as explained above.
Bringing this all together, we can derive the following formula for $R2(Y)$:
\begin{align*}
    R2(Y) & = \int_{0}^{1} \min_{y \in Y} \{ \max(y_1 w, y_2 (1 - w)) \} \, dw \\
                & = \sum_{n=1}^N \left(\color{Orange} u\left(y_1^{(n)}, \left[y_2^{(n)}, y_2^{(n-1)}\right]\right) \color{Black} + \color{NavyBlue} u\left(y_2^{(n)}, \left[y_1^{(n)}, y_1^{(n+1)}\right]\right) \color{Black} \right), \text{where}
\end{align*}
\begin{align*}
    \color{Orange} u\left(y_1^{(n)}, \left[y_2^{(n)}, y_2^{(n-1)}\right]\right) \color{Black} & = \frac 1 2 y_1^{(n)} \left(\left(\frac {y_2^{(n-1)}} {y_1^{(n)} + y_2^{(n-1)}}\right)^2 - \left(\frac {y_2^{(n)}} {y_1^{(n)} + y_2^{(n)}}\right)^2\right), \text{and} \\
    \color{NavyBlue} u\left(y_2^{(n)}, \left[y_1^{(n)}, y_1^{(n+1)}\right]\right) \color{Black} & = \frac 1 2 y_2^{(n)} \left(\left(\frac {y_1^{(n+1)}} {y_1^{(n+1)} + y_2^{(n)}}\right)^2 - \left(\frac {y_1^{(n)}} {y_1^{(n)} + y_2^{(n)}}\right)^2\right).
\end{align*}
The computational complexity of this approach is determined mostly by the condition that we require $Y$ to contain only nondominated solutions and be sorted by $y_1$.
For this, we can first sort an archive of solutions (including dominated points) by $y_1$ and pass over this list once, removing any dominated points and duplicates.
This takes $\mathcal O(N \log N)$ time with standard sorting algorithms.
The computation of the indicator itself is then just a matter of another pass over the sorted list of nondominated points, which requires linear time: All required $w$ and $y$ values can be obtained on the fly based on any given point $y^{(n)}$ and its immediate neighbors.
To summarize, the computation of the exact R2 indicator on an archive of $N$ points in bi-objective space requires a complexity of $\mathcal O(N \log N)$.
This improves upon the $\mathcal O (N |W|)$ complexity of the discretized R2 for precise indicator values and large sets of solutions, i.e., large $|W|$ and $N$ values.

\subsection{Exclusive R2 Contribution} \label{sec:exclusive_r2}

When benchmarking multi-objective optimization algorithms, quality indicators are often not just evaluated after a specific budget of evaluations was spent, but more frequently.
Ideally, the most relevant quality indicator is updated for each change to the nondominated set of points evaluated during a run.
This enables (a) an optimizer to react to large changes in approximation set quality, or (b) a benchmarking framework to record when a target indicator value has been reached for the first time.

For example, the standard protocol of the bi-objective BBOB testbed records the first hitting times for a set of hypervolume indicator values \citep{brockhoff2022using}, using an AVL tree data structure -- i.e., a binary tree whose left and right subtrees differ by at most one level in height \citep{adelson1962avl}  -- to represent the unconstrained nondominated set to represent the unconstrained nondominated set and efficiently updating the hypervolume value for each evaluated point on the test problem.
Each time the approximated Pareto set is updated, the exclusive contributions of removed (added) points is subtracted from (added to) the total indicator value.
Each update is performed in $\mathcal{O} (\log N)$ time, leading to an overall running time of $\mathcal{O} (N \log N)$ to compute the full history of bi-objective R2 indicator values, i.e., the same asymptotic runtime as the individual computation for the final value.

\begin{figure}
    \centering
    \begin{tikzpicture}        
        \coordinate (lu) at (1.3,4.3);
        \coordinate (ym1) at (1.3,3.3);
        \coordinate (yint1) at (2.2,3.3);
        \coordinate (y) at (2.2,2.0);
        \coordinate (yint2) at (3.2,2.0);
        \coordinate (yp1) at (3.2,1.0);
        \coordinate (rl) at (4.2,1.0);
        \coordinate (O) at (0.5,0.5);
        \coordinate (ypm) at (3.2,3.3);
        
        \draw[->,ultra thick] (0,0)--(4,0) node[right, Orange]{$f_1$};
        \draw[->,ultra thick] (0,0)--(0,4) node[above, NavyBlue]{$f_2$};
        \draw[ultra thick, dashed, LightGray] (lu) -- (ym1);
        \draw[ultra thick, LightGray] (ym1) -- (yint1);
        \draw[ultra thick, Orange] (yint1) -- (y);
        \draw[ultra thick, NavyBlue] (y) -- (yint2);
        \draw[ultra thick, LightGray] (yint2) -- (yp1);
        \draw[ultra thick, dashed, LightGray] (yp1) -- (rl);
        
        \draw[ultra thick, dashed, NavyBlue] (yint1) -- (ypm);
        \draw[ultra thick, dashed, Orange] (ypm) -- (yint2);

        \draw (ym1) node {\huge +} node[below left]{$y^{(n-1)}$};
        % \draw (yint1) node {\huge +};
        \draw (y) node {\huge +} node[below left]{$y^{(n)}$};
        % \draw (yint2) node {\huge +};
        \draw (yp1) node {\huge +} node[below left]{$y^{(n+1)}$};
        \draw[gray] (O) node {\huge +};
    \end{tikzpicture}
    
    % \caption{Values used while computing the exclusive R2 contribution of $y^{(n)}$. $y^{(n \pm)}$~contains the respective maximum values from $y^{(n-1)}$ and $y^{(n+1)}$, and identifies the knee in the Pareto front that is introduced when $y^{(n)}$ is removed. For each $y$, there is a corresponding weight vector $w$.}
    \caption{Illustration of the effect of removing $y^{(n)}$ from $Y$ to determine the exclusive contribution $\Delta R2(y^{(n)},Y)$. Both PF line segments with constant $\color{Orange} f_1$ and $\color{NavyBlue} f_2$ values are shifted upwards to the dashed lines.}
    \label{fig:r2-contr}
\end{figure}

To implement an iterative version of the R2 indicator, we most importantly need to be able to compute the exclusive contribution $\Delta R2(y^{(n)}, Y)$ for a solution $y^{(n)} \in Y$ to the overall value.
We consider the effect of removing $y^{(n)}$ from $Y$, cf.~\Cref{fig:r2-contr}:
\begin{align*}
    \Delta R2(y^{(n)}, Y) & = R2(Y) - R2(Y \setminus \{y^{(n)}\}) \\
    & = \color{Orange} u\left(y_1^{(n)},[y_2^{(n)},y_2^{(n-1)}]\right) - u\left(y_1^{(n+1)},[y_2^{(n)},y_2^{(n-1)}]\right) \color{Black} \\
    & + \color{NavyBlue} u\left(y_2^{(n)},[y_1^{(n)},y_1^{(n+1)}]\right) - u\left(y_2^{(n-1)},[y_1^{(n)},y_1^{(n+1)}]\right) \color{Black}. % \\
    % & = \dots?
     % &= R2(y^{(n)}) - \int_{w^{(n-)}_1}^{w^{(n+)}_1} \min{(y^{(n)}_1 w, y^{(n)}_2 (1 - w))} \, dw \\
     % &= R2(y^{(n)}) - \left(\int_{w^{(n-)}_1}^{w^{(n \pm)}_1} y^{(n)}_1 w \, dw + \int_{w^{(n \pm)}_1}^{w^{(n+)}_1} y^{(n)}_2 (1 - w) \, dw \right) \\
     % &= R2(y^{(n)}) - \left(\left[\frac 1 2 y^{(n)}_1 w^2\right]_{w^{(n-)}_1}^{w^{(n \pm)}_1} + \left[- \frac 1 2 y^{(n)}_2 (1 - w)^2\right]_{w^{(n \pm)}_1}^{w^{(n+)}_1} \right) \\
     % &= R2(y^{(n)}) - \left(\frac 1 2 y^{(n)}_1 \left((w^{(n \pm)}_1)^2 - (w^{(n-)}_1)^2\right) + \frac 1 2 y^{(n)}_2 \left((1 - w^{(n \pm)}_1)^2 - (1 - w^{(n+)}_1)^2\right) \right)
\end{align*}
%$\Delta R2(y^{(n)}, Y)$ is negative if there is more than one point in $Y$, as the corresponding Tchebycheff aggregation functions have smaller values along the Pareto front more closer to the ideal point, i.e., $R2(Y) < R2(Y \setminus \{y\})$.
Note that a Pareto front $A$ that (weakly) dominantes another Pareto front $B$ will be closer to the ideal point and thus
reduce $R2$ as the corresponding Tchebycheff aggregation functions have smaller values. This is in contrast to the HV indicator, where better fronts yield in an increased value. 
In consequence, except for the special case $Y = \{y\}$, $R2(Y) < R2(Y \setminus \{y\})$ and thus $\Delta R2(y^{(n)}, Y)$ will be negative.
Further, this defines the region of the exclusive contribution of each solution in $Y$, and so re-establishes Pareto compliance for the R2 indicator as it was described in its original publication \citep{hansen1998}.

\subsection{Incremental Bi-objective R2 Indicator Computation} \label{sec:incremental_r2}

Using the exclusive contribution $\Delta R2(y^{(n)}, Y)$, the procedure \textsc{R2IndicatorUpdate} in \Cref{alg:update} demonstrates how an update of the nondominated set $Y$ and the corresponding R2 indicator value can be performed.
To start with, we check whether $Y$ is empty, in which case we add $y$ to it and compute the R2 indicator value for this individual point.
Otherwise, if $y$ is not dominated by some point in $Y$, we remove all points dominated by $y$ from $Y$, and finally add $y$ to $Y$.
Using the $\Delta R2$ contribution, we update the indicator correspondingly for these removal and addition operations.

% To start with, we efficiently find the solution $y^{(i)}$ with the lowest $f_1$ value higher than or equal to $y$ in $Y$.
% Then we can iteratively remove all solutions dominated by $y$, subtracting their exclusive R2 contributions to $Y'$ from $R2'$, respectively.
% Finally, if $y$ is now non-dominated to the next solution, we add it to the next position in $Y$ and add its corresponding R2 contribution.

% If $y$ is dominated by or equals some solution in $Y$, nothing happens, as both modifications to $Y'$ and $R2'$ are not performed.
% If $y$ dominates some solutions in $Y$, they are iteratively identified and removed, before inserting $y$ into the sorted set.
% If $y$ is mutually nondominated with every solution in $Y$, no solutions are removed, and only the latter modification - addition to the set - is performed.

\begin{algorithm}[!t]
    \caption{Bi-objective R2 Indicator Update}\label{alg:update}
    \begin{algorithmic}[1]
        
        \Input {Nondominated solutions $Y$, $R2(Y)$, new solution $y \in \mathbb R^2$}
        \EndInput
        
        \Output{Updated nondominated solutions $Y'$, $R2(Y')$}
        \EndOutput
        
        \Procedure{R2IndicatorUpdate}{$Y$, $R2(Y), y$}
            
            \State $Y' \gets Y$
            \State $R2' \gets R2(Y)$

            \If {$Y = \varnothing$}
                \State $R2' \gets R2(y)$
                \State $Y' \gets \{y\}$
            \ElsIf {$\nexists y' \in Y:  y' \preceq y$} \Comment{Proceed if $y$ is not dominated by $Y$}
                \For{$y' \in \{y^\star \in Y \mid y \preceq y^\star\}$} \Comment{Remove all solutions dominated by $y$}
                    \State $R2' \gets R2' - \Delta R2(y', Y)$
                    \State $Y' \gets Y' \setminus \{y'\}$
                \EndFor
                \State{$Y' \gets Y' \cup \{y\}$} \Comment{Add $y$ to $Y'$}
                \State{$R2' \gets R2' + \Delta R2(y, Y)$}
            \EndIf

            % \If{$Y = \varnothing$}
            %     \State $R2' = R2(y)$
            %     \State $Y' \gets \{y\}$
            % \Else
            %     \State $i \gets$ \Call{FindLower}{$Y$, $y$}
                
            %     \While{$y \prec y^{(i)}$} \Comment{Remove all solutions dominated by $y$}
            %         \State $R2' \gets R2' - \Delta R2(y^{(i)}, Y)$
            %         \State $Y' \gets Y' \setminus \{y^{(i)}\}$ \Comment{$y^{(i)}$ now points to the next solution in $Y$}
            %     \EndWhile
    
            %     \If{$y^{(i)} \npreceq y$ \textbf{and} $y \npreceq y^{(i)}$} \Comment{Add $y$ to $Y$ if nondominated}
            %         \State $Y' \gets Y' \cup \{y\}$
            %         \State $R2' \gets R2' + \Delta R2(y, Y)$
            %     \EndIf
                
            % \EndIf
                            
            \State \textbf{return} $(Y', R2')$
        \EndProcedure
    \end{algorithmic}
\end{algorithm}

Storing the nondominated set $Y$ with a suitable data structure for a sorted list (e.g., AVL trees), we can efficiently determine whether $y$ is dominated by some point from $Y$ in $\mathcal O(\log N)$ time.
Addition and deletion from the nondominated set can equally be performed in $\mathcal O(\log N)$ for each operation, and as each point can only be added or deleted once, we get an overall runtime of $\mathcal O(N \log N)$ to perform $N$ successive \textsc{R2IndicatorUpdate} operations to compute the complete history of R2 indicator values across $N$ evaluated solutions.

\section{Properties of the Exact R2 Indicator} \label{sec:experimental}

In this section, we present some empirical and theoretical results on our proposed exact R2 indicator.
We start by demonstrating the approximation behavior of the discrete R2 in relation to the exact computation described in the previous sections.
Then, we will calculate the optimal R2 indicator values for simple front shapes and demonstrate the convergence behavior of the indicator w.r.t. increasingly large nondominated sets.

In \cite{schaepermeier2024reinvestigating}, we implemented the exact R2 indicator in the statistical software \texttt{R}.
We subsequently implemented incremental variants of the R2 and hypervolume indicators in Python due to the easy access to an efficient sorted list in the \texttt{sortedcontainers} library \citep{jenks2025sorted}.
For computations of the discrete R2 indicator, we utilize the \texttt{unary\_r2\_indicator} function from the \texttt{emoa} \texttt{R} package \citep{mersmann2023emoa}.
The scripts to reproduce these experimental results are published at \url{https://github.com/schaepermeier/r2-revisited}.

\subsection{Comparison of Discrete and Exact R2 Values}

\label{sec:discrete}

\begin{figure}[tb]
    \centering
    \includegraphics[width=.48\linewidth]{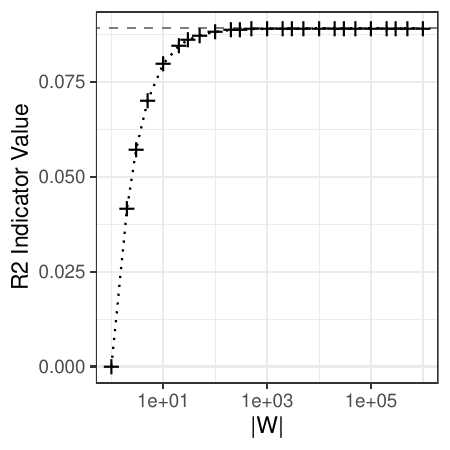}
    \hfill
    \includegraphics[width=.48\linewidth]{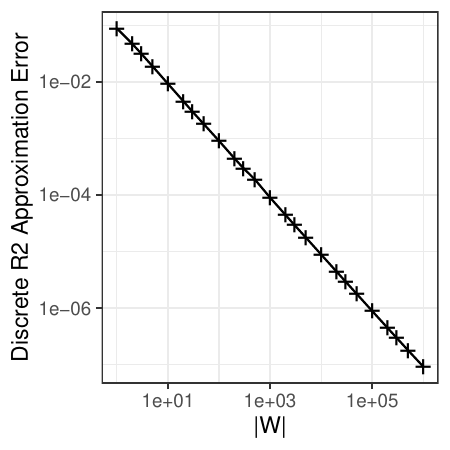}
    \caption{Comparison of the approximated R2 indicator values using $|W|$ weight vectors and the exact R2 indicator value computed by our methodology. Left: The exact value is shown by the dashed line. Right: The corresponding approximation error of the discretized approach. The solution set is given by $1,\!001$ equidistant solutions evaluated on the Pareto set of a bi-sphere problem with centers $(-0.5, 0)$ and $(0.5, 0)$.}
    \label{fig:r2-approx}
\end{figure}

We start by comparing the exact R2 indicator as computed using our method (see \Cref{sec:method}) with the discretized version found throughout the literature.
As a basis for the comparison, we sample equidistant solutions on the Pareto set of a bi-sphere problem $F(x) = \left(\sum (x - c_1)^2, \sum (x - c_2)^2\right)$ with centers $c_1 = (-0.5, 0)$ and $c_2 = (0.5, 0)$ using a step size of $0.001$, resulting in $1,\!001$ points.
We evaluated the discretized R2 with uniformly distributed weights, and the number of weights~($|W|$) ranging from one to one million.
The value of the discretized indicator as well as the approximation error are visualized in \Cref{fig:r2-approx}.

We can see that, in this example, the R2 indicator seems sufficiently well approximated at around $1,\!000$ weights.
Still, more weights yield a more accurate approximation: Empirically, there seems to be an exponential relationship between the number of weights chosen and the approximation error.
This reflects analyses by \cite{brockhoff2012properties} on the behavior of the discrete R2 with an increasing number of weights, albeit missing the exact R2 values for comparison.

\subsection{Exact R2 Indicator Values}

For the exact R2 indicator, we can provide the optimal indicator values for simple solution sets and corresponding test functions.
This contributes to a better understanding of the R2 indicator values, as a geometric interpretation like with the hypervolume indicator is not possible.

\subsubsection{Nadir and Ideal Points}

Assuming we have normalized the objective space ($(0,0)$ being the ideal point and $(1,1)$ being the nadir point), we can compute the R2 indicator for the worst possible solution within the region of interest $[0,1] \times [0,1]$ by inserting the nadir point $(1,1)$ into the equation. 
According to the previous results, we can compute this value as follows: 
\begin{align*}
    % R2(\{(1,1)\}) &= \frac 1 2 y_2 (1 - (1 - w^*)^2) + \frac 1 2 y_1 (1 - (w^*)^2) \\
    % & = 0.5 \cdot 1 (1 - (1 - 0.5)^2) + 0.5 \cdot 1 (1 - 0.5^2) \\
    % &= 0.5 \cdot 0.75 + 0.5 \cdot 0.75 \ \mathbf{= 0.75}.
    R2(\{(1,1)\}) &= \frac 1 2 y_1 \left(1 - \left(\frac {y_2} {y_1 + y_2}\right)^2\right) + \frac 1 2 y_2 \left(1 - \left(\frac {y_1} {y_1 + y_2}\right)^2\right) \\
    & = 0.5 \cdot (1 - 0.25) + 0.5 \cdot (1 - 0.25) \ \mathbf{= 0.75}
\end{align*}
This result is independent of the particular problem or PF, as it is only dependent on the normalized nadir and ideal points.
Analogously, we can derive the R2 value for the ideal point as $R2(\{(0,0)\}) \ \mathbf{= 0}$, as all utilities equal zero at the ideal point.
For comparison, the HV w.r.t. the nadir point as the reference point is $HV(\{(1, 1)\}) = 0$, while the HV of the ideal point is $HV(\{(0, 0)\}) = 1$ in this situation. 

\subsubsection{Linear Front}

Let us now consider a linear PF where $y_2 = 1 - y_1$ and $y_1 \in [0,1]$ with ideal point $(0,0)$.
When computing R2, for each weight vector $w$, we can find the optimal solution $y$ on the PF, which we can derive as follows:
\begin{align*}
    w_1 y_1 = w_2 y_2 & \Rightarrow w_1 y_1 = (1 - w_1) (1 - y_1) \Leftrightarrow y_1 = 1 - w_1.
\end{align*}
Note that when we use this, it does not matter which of the objectives we consider in the R2 computation, as they will always yield the same utility value.
Integrating over the weights for this set $Y_\texttt{lin}$, we get:
\begin{align*}
    R2(Y_\texttt{lin}) &= \int_0^1 \min_{y \in Y_\texttt{lin}} \{\max(w y_1, (1 - w) y_2)\} \, dw \\
          &= \int_0^1 w y_1 \, dw = \int_0^1 w (1 - w) \, dw = \int_0^1 (w - w^2) \, dw \\
          &= \left[\frac 1 2 w^2 - \frac 1 3 w^3 \right]^1_0 \\
          &= \frac 1 2 - \frac 1 3 \ \mathbf{= \frac 1 6 \approx 0.1667}.
\end{align*}

Based on this result, we can derive the optimal R2 indicator value for any problem with a linear PF.
For example, the well-known DTLZ1 problem \citep{deb2005scalable} possesses a linear PF with ideal point $(0, 0)$ and nadir point $(0.5, 0.5)$, which results in an optimal R2 value of $\frac 1 {12} \approx 0.0833$.

\subsubsection{Convex and Concave Fronts}

\begin{figure}[t]
    \centering
    \includegraphics[width=.32\linewidth]{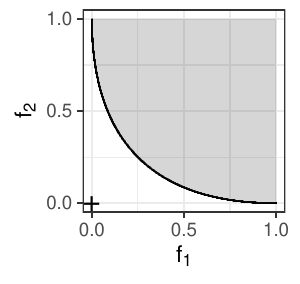}
    \includegraphics[width=.32\linewidth]{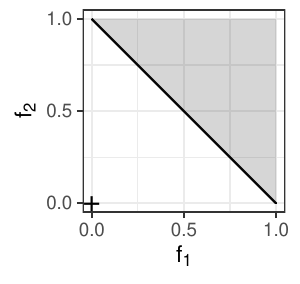}
    \includegraphics[width=.32\linewidth]{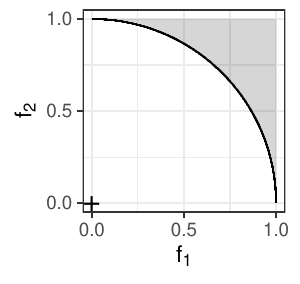}
    \caption{Schematic illustration of a convex ($0 < R2 < \frac 1 6$), linear ($R2 = \frac 1 6$), and concave ($\frac 1 6 < R2 < \frac 1 4 $) PF, respectively. The ideal point at the origin is denoted by \texttt{+}, and the gray area indicates the dominated area.}
    \label{fig:schematic-fronts}
\end{figure}

Additionally, we can derive value ranges for general concave and convex PFs: A concave front $Y_\texttt{conc}$ will always achieve worse utility values than a linear function, and a convex front $Y_\texttt{conv}$ will always have better utility.
Again considering the normalized objective space, a convex front $Y_\texttt{conv}$ will always fall between the R2 values of the ideal point and the linear front, i.e., $0 < R2(Y_\texttt{conv}) < \frac 1 6$.
Analogously, a general concave front has an ideal R2 value between the value of the linear front and a front made up of only the extreme points $\{(0,1), (1,0)\}$.
This extreme concave front has an R2 value of $u(1, [0, 1]) + u(1, [0, 1]) = 2 \cdot 0.5 \cdot (0.5^2 - 0^2) = 0.25$, so $\frac 1 6 < R2(Y_\texttt{conc}) < \frac 1 4$.
All cases are illustrated in \Cref{fig:schematic-fronts}.
If none of the discussed conditions apply, the ideal R2 indicator value of the normalized objective space may lie anywhere between $0$ and $0.25$.

As shown above, exact R2 indicator values can be derived for certain analytical PF shapes.
Following the same pattern, that is, resolving $w_1 y_1 = (1 - w_1) y_2$ and integrating the utility function, we can compute the exact indicator values also for more complex PF shapes, albeit in a less straightforward manner.
We limit ourselves to reporting the results for simple quadratic PF functions, which correspond to the PFs in \Cref{fig:schematic-fronts}:
\begin{itemize}
    \item Convex PF with $y_2 = (1 - \sqrt{y_1})^2$: $\frac {3 \pi - 8} {16} \approx 0.0890$
    \item Concave PF with $y_2 = \sqrt{1 - y_1^2}$: $\frac 1 8 \left(3 \sqrt 2 \sinh^{-1}{(1)} - 2\right) \approx 0.2174$
\end{itemize}
The result for the convex PF applies, e.g., for the classical bi-sphere problem, while the concave PF corresponds to DTLZ2-6 \citep{deb2005scalable}.

To summarize, we can compute exact R2 indicator values for different simple front shapes and individual points.
While these values do not seem to have an intuitive (geometric) interpretation, they are rather given meaning by the expected utility to a decision-maker.

\subsection{Iterative Indicator Comparison} \label{sec:iterative_comparison}

\begin{figure}[t]
    \centering
    \begin{subfigure}{0.49\linewidth}
        \includegraphics[width=\linewidth]{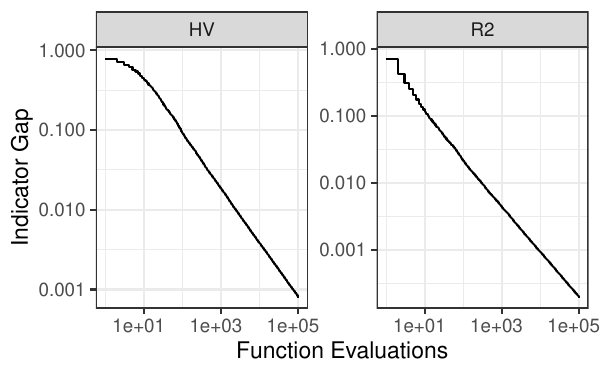}
        \caption{Random uniform sampling in $[-1,1]^2$.}
    \end{subfigure}
    \hfill
    \begin{subfigure}{0.49\linewidth}
        \includegraphics[width=\linewidth]{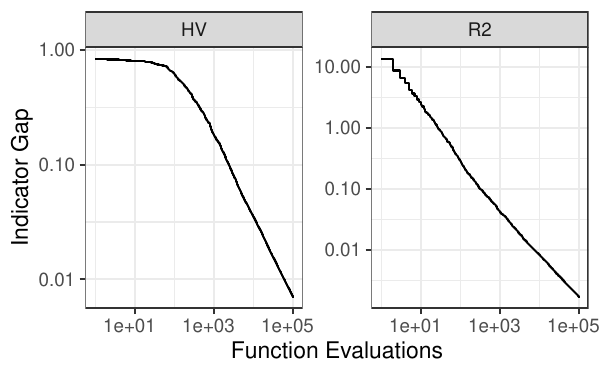}
        \caption{Random uniform sampling in $[-5,5]^2$.}
    \end{subfigure}
    \caption{Comparison of the hypervolume (HV) and R2 indicator gaps using the nadir and ideal point, respectively, on a bi-sphere problem with centers $(-0.5, 0)$ and $(0.5, 0)$. As the R2 indicator always considers all solutions dominated by the ideal point, it does not produce an initial plateau.}
    \label{fig:iterative-indicator}
\end{figure}

Finally, we compare the gaps of the hypervolume and R2 indicators, respectively, when evaluating a large amount of randomly sampled solutions, cf.~\Cref{fig:iterative-indicator}.
We use the same bi-sphere problem from above (see \Cref{sec:discrete}) with centers $(-0.5,0)$ and $(0.5,0)$, and sample $10^5$ solutions from $[-5,5]^2$ and $[-1,1]^2$ uniformly at random, respectively.
The nadir ($1,1$) and ideal points $(0,0)$ are used as the respective reference points, and average indicator gaps after $100$ repetitions are shown.

Here, we can see how not considering all solutions in the indicator computation affects the HV indicator: It produces a noticeable plateau with few sampled solutions, which is all the more pronounced when the sampling area is larger.
This emphasizes the practicality of the R2 indicator for benchmarking scenarios where an ideal point, but not necessarily a good anti-optimal reference point, can be generated automatically.

\section{Conclusion} \label{sec:conclusion}

In this paper, we introduce a procedure to compute the exact R2 indicator value for a given solution set.
In contrast to its widely-known and commonly used discrete counterpart, the exact R2 indicator is not just weakly Pareto-compliant, but a proper Pareto-compliant indicator.
This is achieved by foregoing the discretization of the distribution of utility functions usually performed in the calculation of the indicator and taking a continuous, uniform distribution of Tchebycheff utility functions as the basis.
Pareto compliance of the R2 indicator with this utility distribution was already described when it was introduced by \cite{hansen1998}, but missing instructions for its exact calculation.
Further, we show how this indicator is implemented efficiently with a running time of $\mathcal O(N \log N)$ for bi-objective solution sets of size $N$.
Further, we derive the exclusive contribution of an individual solution to the indicator value, which we then use in an incremental update procedure to determine the full history of indicator values within the same runtime of $\mathcal O(N \log N)$.
This positions the exact R2 indicator as a promising Pareto-compliant alternative to the hypervolume indicator, especially when a utopian rather than an anti-optimal reference point is naturally available, making it the more natural choice for benchmarking multi-objective optimizers in common problem settings.
We demonstrate the approximation behavior of the commonly used, discretized R2 indicator in comparison with the exact computation, and provide optimal R2 indicator values for the ideal and nadir points, as well as a linear front in normalized objective space.
We then use these values to derive R2 indicator value ranges for general convex and concave PFs as well.

The availability of the exact, and thus Pareto-compliant, R2 indicator along with the proposed approach to compute it efficiently offers multiple directions for further theoretical and empirical research.
% So far, we have only demonstrated how the R2 indicator is computed for bi-objective problems.
% A natural further research direction pertains to the computation of the indicator for more than two objectives.
% We do not expect that the R2 indicator will provide a runtime advantage over the hypervolume indicator, however, schemes to approximate it are until now the state of the art in its computation and therefore very accessible compared to HV approximations.
%
From a theoretical point of view, we see potential in analyzing the approximation quality of the discretized R2 depending on the number of weights used.
This could be particularly interesting for giving quality guarantees for R2 in higher-dimensional objective spaces, where exact indicator computations may become computationally intractable.
A deeper theoretically supported analysis of its properties, along the lines of \cite{brockhoff2012properties}, would also be interesting.
Further, the effects of different utility functions as well as the integration of (decision-maker) preferences can be examined, as both directions have been studied in detail for the discretized R2 indicator \citep{wagner2013preference}.
Finally, connections between the R2 indicator and the integrated preference functional \citep{carlyle2003quantitative,bozkurt2010quantitative}, a parallel development of an indicator very similar to R2 in the operations research community, should be further investigated, and could yield improvements in understanding the R2 indicator's properties and computation.

Differences and similarities in the preferred distributions between the exact R2 indicator and the hypervolume indicator when applying them in a benchmarking context could present a further promising research direction.
Likewise, integrating the exact R2 in optimization heuristics, similar to R2-EMOA \citep{trautmann2013r2}, may be the subject of future studies.

% ===== Bibliography =====

\small

\bibliographystyle{apalike}
\bibliography{library}

\end{document}